\documentclass[12pt]{amsart}
\usepackage{amsmath, amssymb, epsfig}

\newlength{\algorithmwidth}
\theoremstyle{plain}
\newtheorem{theorem}{Theorem}[section]
\newtheorem{proposition}[theorem]{Proposition}
\newtheorem{corollary}[theorem]{Corollary}
\newtheorem{lemma}[theorem]{Lemma}

\theoremstyle{definition}
\newtheorem{definition}[theorem]{Definition}

\theoremstyle{remark}

\numberwithin{equation}{section}

\DeclareMathOperator*{\supp}{supp}

\DeclareMathOperator*{\range}{range}

\DeclareMathOperator*{\argmin}{arg min}

\def \R {\mathbb{R}}

\def \e {\varepsilon}

\def \d {\delta}

\def \L {\Lambda}

\def \< {\langle}
\def \> {\rangle}
\def \^ {\widehat}

\def \supp {{\rm supp}}

\begin{document}
\title[]{Signal recovery from incomplete and inaccurate measurements 
  via Regularized Orthogonal Matching Pursuit}

\author{Deanna Needell
  \and Roman Vershynin}

\thanks{Partially supported by the Alfred P.~Sloan Foundation
  and by NSF DMS grant 0652617}
  
\date{October 5, 2007}

\address{Department of Mathematics,
   University of California,
   Davis, CA 95616, USA}
\email{\{dneedell,vershynin\}@math.ucdavis.edu}

\begin{abstract}
  We demonstrate a simple greedy algorithm that can reliably recover a vector 
  $v \in \R^d$ from incomplete and inaccurate measurements 
  $x = \Phi v + e$. Here $\Phi$ is a $N \times d$ measurement matrix with $N \ll d$,
  and $e$ is an error vector. Our algorithm, Regularized Orthogonal Matching 
  Pursuit (ROMP), seeks to close the gap between two major approaches to sparse 
  recovery. It combines the speed and ease of implementation of the greedy methods 
  with the strong guarantees of the convex programming methods.

  For any measurement matrix $\Phi$ that satisfies a Uniform Uncertainty Principle, 
  ROMP recovers a signal $v$ with $O(n)$ nonzeros from its inaccurate measurements $x$ 
  in at most $n$ iterations, where each iteration amounts to solving a Least Squares Problem.
  The noise level of the recovery is proportional to $\sqrt{\log n} \|e\|_2$. 
  In particular, if the error term $e$ vanishes the reconstruction is exact.
  
  This stability result extends naturally to the very accurate recovery of 
  approximately sparse signals.
\end{abstract}
\subjclass{68W20, 65T50, 41A46}
\maketitle

\section{Introduction}

\subsection{Exact recovery by convex programming}

The recent massive work in the area of Compressed Sensing, surveyed in \cite{C},
rigorously demonstrated that one can algorithmically 
recover sparse (and, more generally, compressible) signals from incomplete observations.
The simplest model is a $d$-dimensional signal $v$ with a small number of nonzeros:
$$
v\in\R^d, \qquad |\supp(v)|\leq n \ll d.
$$ 
Such signals are called $n$-sparse.
We collect $N\ll d$ nonadaptive linear measurements of $v$,
given as $x = \Phi v$ where $\Phi$ is some $N$ by $d$ measurement matrix. 
We then wish to efficiently recover the signal $v$ from its measurements $x$. 

A necessary and sufficient condition for exact recovery is that the map $\Phi$ 
be one-to-one on the set of $n$-sparse vectors. 
Cand\`es and Tao \cite{CT decoding} proved that under a stronger (quantitative) 
condition, the sparse recovery problem is equivalent to a convex program
\begin{equation}                    \label{L1}
  \min \|u\|_1 \qquad \text{subject to} \qquad \Phi u = x
\end{equation}
and therefore is computationally tractable. This condition is that the map
$\Phi$ is an almost isometry on the set of $O(n)$-sparse vectors:

\begin{definition}[Restricted Isometry Condition] 
  A measurement matrix $\Phi$ satisfies the 
  {\em Restricted Isometry Condition} (RIC)
  with parameters $(m, \e)$ for $\e \in (0,1)$
  if we have
  $$
  (1-\e)\|v\|_2 \leq \|\Phi v\|_2 \leq (1+\e)\|v\|_2
  \qquad \text{for all $m$-sparse vectors}.
  $$
\end{definition}

Under the Restricted Isometry Condition with parameters $(3n, 0.2)$, 
the convex program \eqref{L1} exactly recovers an $n$-sparse signal $v$
from its measurements $x$ \cite{CT decoding}.

The Restricted Isometry Condition can be viewed as an abstract 
form of the Uniform Uncertainty Principle of harmonic analysis
(\cite{CT}, see also \cite{CRT} and \cite{LV}).
Many natural ensembles of random matrices, such as partial Fourier, 
Bernoulli and Gaussian, satisfy the Restricted Isometry condition with 
parameters $n \ge 1$, $\e \in (0,1/2)$ provided that
$$
N = n \e^{-O(1)} \log^{O(1)} d;
$$
see e.g. Section 2 of \cite{NV} and the references therein.
Therefore, a computationally tractable exact recovery of sparse signals is 
possible with the number of measurements $N$ roughly proportional to the 
sparsity level $n$, which is usually much smaller than the dimension $d$.

\subsection{Exact recovery by greedy algorithms}

An important alternative to convex programming is greedy algorithms, which have
roots in Approximation Theory. A greedy algorithm computes the support of $v$ iteratively, 
at each step finding one or more new elements (based on some ``greedy'' rule) 
and subtracting their contribution from the measurement vector $x$.
The greedy rules vary. The simplest rule is to pick a coordinate of $\Phi^* x$ of 
the biggest magnitude; this defines the well known greedy algorithm called 
Orthogonal Matching Pursuit (OMP), known otherwise as 
Orthogonal Greedy Algorithm (OGA) \cite{TG}.

Greedy methods are usually fast and easy to implement, which makes them popular 
with practitioners. For example, OMP needs just $n$ iterations to find the
support of an $n$-sparse signal $v$, and each iteration amounts to solving one least-squares 
problem; so its running time is always polynomial in $n$, $N$ and $d$.
In contrast, no known bounds are known on the running time of \eqref{L1}
as a linear program. Future work on customization of convex programming
solvers for sparse recovery problems may change this picture, of course.
For more discussion, see \cite{TG} and \cite{NV}.

A variant of OMP was recently found in \cite{NV} that has guarantees essentially 
as strong as those of convex programming methods.\footnote{OMP itself does not have such strong guarantees, see \cite{R}.}
This greedy algorithm is called Regularized Orthogonal Matching Pursuit (ROMP); 
we state it in Section~\ref{s: stable} below.
Under the Restricted Isometry Condition with parameters $(2n, 0.03 / \sqrt{\log n})$, ROMP exactly recovers an $n$-sparse signal $v$ from its measurements $x$.  

Summarizing, 
{\em the Uniform Uncertainty Principle is a guarantee for efficient sparse recovery;
one can provably use either convex programming methods \eqref{L1}
or greedy algorithms (ROMP).}

\subsection{Stable recovery by convex programming and greedy algorithms}  \label{s: stable}

A more realistic scenario is where the measurements are inaccurate (e.g. contaminated by noise) and the signals are not exactly sparse. 
In most situations that arise in practice, one cannot hope to know the measurement vector
$x = \Phi v$ with arbitrary precision. Instead, it is perturbed by a small
error vector: $x = \Phi v + e$. Here the vector $e$ has unknown coordinates as well as unknown magnitude, and it needs not be sparse (as all coordinates may be affected by the noise). For a recovery algorithm to be stable, it should be able to approximately 
recover the original signal $v$ from these perturbed measurements. 

The stability of convex optimization algorithms for sparse recovery was studied
in \cite{D}, \cite{T}, \cite{DET}, \cite{CRT-stability}. 
Assuming that one knows a bound on the magnitude 
of the error, $\|e\| \le \d$, it was shown in \cite{CRT-stability} that 
the solution $\hat{v}$ of the convex program 
\begin{equation}                    \label{L1-regularized}
  \min \|u\|_1 \qquad \text{subject to} \qquad \|\Phi u - x\|_2 \le \d
\end{equation}
is a good approximation to the unknown signal: 
$\|v - \hat{v}\|_2 \le C \d$.

In contrast, the stability of greedy algorithms for sparse recovery has not been
well understood. Numerical evidence \cite{DET} suggests that OMP should be less 
stable than the convex program \eqref{L1-regularized}, but no theoretical 
results have been known in either the positive or negative direction.
The present paper seeks to remedy this situation. 

We prove that {\em ROMP is as stable as the convex program \eqref{L1-regularized}.}
This result essentially closes a gap between convex programming and greedy
approaches to sparse recovery.

\bigskip

\textsc{Regularized Orthogonal Matching Pursuit (ROMP)}

\nopagebreak

\fbox{\parbox{\algorithmwidth}{
  \textsc{Input:} Measurement vector $x \in \R^N$ and sparsity level $n$
  
  \textsc{Output:} Index set $I \subset \{1,\ldots,d\}$, reconstructed vector $\hat{v}=y$

  \begin{description}
    \item[Initialize] Let the index set $I = \emptyset$ and the residual $r = x$.\\
      Repeat the following steps $n$ times or until $|I|\geq 2n$:
      \begin{description}
    \item[Identify] Choose a set $J$ of the $n$ biggest nonzero coordinates in magnitude 
      of the observation vector $u = \Phi^*r$, or all of its nonzero coordinates, 
      whichever set is smaller.
    \item[Regularize] Among all subsets $J_0 \subset J$ with comparable coordinates:
      $$
      |u(i)| \leq 2|u(j)| \quad \text{for all } i,j \in J_0,
      $$
      choose $J_0$ with the maximal energy $\|u|_{J_0}\|_2$.
    \item[Update] Add the set $J_0$ to the index set: $I \leftarrow I \cup J_0$, 
      and update the residual:
      $$
      y = \argmin_{z \in \R^I} \|x - \Phi z\|_2; \qquad r = x - \Phi y.
      $$
      \end{description}
      
  \end{description}
 }}

\bigskip


\begin{theorem}[Stability under measurement perturbations]\label{T:stability}
  Assume a measurement matrix $\Phi$ satisfies the Restricted Isometry Condition 
  with parameters $(4n, \e)$ for $\e = 0.01 / \sqrt{\log n}$. 
  Let $v$ be an $n$-sparse vector in $\R^d$. 
  Suppose that the measurement vector $\Phi v$ becomes corrupted, so we consider
  $x = \Phi v + e$ where $e$ is some error vector. 
  Then ROMP produces a good approximation to $v$:
  $$
  \|v - \hat{v}\|_2 \leq 104 \sqrt{\log n}\|e\|_2.
  $$  
\end{theorem}

Note that in the noiseless situation ($e = 0$) the reconstruction is exact:
$\hat{v} = v$. This case of Theorem~\ref{T:stability} was proved in \cite{NV}.

Our stability result extends naturally to the even more realistic scenario 
where the signals are only approximately sparse. Here and henceforth, denote by $f_m$ the vector of the $m$ biggest coefficients in absolute value of $f$. 

\begin{corollary}[Stability of ROMP under signal perturbations]\label{T:stabsig}
  Assume a measurement matrix $\Phi$ satisfies the Restricted Isometry Condition 
  with parameters $(8n, \e)$ for $\e = 0.01 / \sqrt{\log n}$. 
  Consider an arbitrary vector $v$ in $\R^d$.
 Suppose that the measurement vector $\Phi v$ becomes corrupted, 
  so we consider $x = \Phi v + e$ where $e$ is some error vector. 
  Then ROMP produces a good approximation to $v_{2n}$:
  \begin{equation}\label{boundsig}
    \|\hat{v} - v_{2n}\|_2 
    \leq 159 \sqrt{\log 2n} \Big( \|e\|_2 + \frac{\|v-v_n\|_1}{\sqrt{n}} \Big).
  \end{equation}  
\end{corollary}

\begin{remarks}
{\bf 1. } The term $v_{2n}$ in the corollary can be replaced by $v_{(1+\delta)n}$ for any $\delta > 0$. This change will only affect the constant terms in the corollary.

  {\bf 2. } By applying Corollary~\ref{T:stabsig} 
    to the largest $2n$ coordinates of $v$ and using Lemma~\ref{L:ve} below, 
    we also have the error bound for the entire vector $v$:
    \begin{equation}\label{vbound}
      \|\hat{v} - v\|_2 
      \leq 160 \sqrt{\log 2n} \Big( \|e\|_2 + \frac{\|v-v_n\|_1}{\sqrt{n}} \Big).
    \end{equation}

  {\bf 3. } For the convex programming method \eqref{L1-regularized}, 
    the stability bound \eqref{vbound} was proved in \cite{CRT-stability}, 
    and even without the logarithmic factor.
    We conjecture that this factor is also not needed in our results for ROMP.

  {\bf 4. } Unlike the convex program \eqref{L1-regularized}, ROMP succeeds
    with absolutely no prior knowledge about the error $e$; 
    its magnitude can be arbitrary. In the terminology of \cite{DET}, 
    the convex programming approach needs to be ``noise-aware'' while ROMP needs not. 

  {\bf 5. } One can use ROMP to approximately compute a $2n$-sparse vector that is close to {\em the best $2n$-term approximation}
    $v_{2n}$ of an arbitrary signal $v$. To this end, one just needs to retain the $2n$ biggest coordinates of $\hat{v}$.
    Indeed, Corollary~\ref{C:napprox} below shows that the best $2n$-term approximations
    of the original and the reconstructed signals are close:
    $$
    \|v_{2n} - \hat{v}_{2n}\|_2 \leq 477 \sqrt{\log 2n}\Big( \|e\|_2 + \frac{\|v-v_{n}\|_1}{\sqrt{n}}\Big).
    $$

  {\bf 6. } An important special case of Corollary~\ref{T:stabsig} is for the class of 
    compressible vectors, which is a common model in signal processing, see \cite{CT}, \cite{Do}. 
    Suppose $v$ is a compressible vector in the sense that its coefficients 
    obey a power law: for some $p > 1$, the $k$-th largest coefficient in magnitude of $v$
    is bounded by $C_p k^{-p}$.
    Then \eqref{vbound} yields the following bound on the reconstructed signal:
    \begin{equation}\label{boundcomp}
      \|v - \hat{v}\|_2 \le C'_p \frac{\sqrt{\log n}}{n^{p-1/2}} + C''\sqrt{\log n}\|e\|_2.
    \end{equation} 
    As observed in \cite{CRT-stability}, this bound is optimal (within the logarithmic
    factor); no algorithm can perform fundamentally better.
\end{remarks}

The rest of the paper is organized as follows. 
In Section~\ref{s: proof}, we prove our main result, Theorem~\ref{T:stability}.
In Section~\ref{s: consequences}, we deduce the extension for 
approximately sparse signals, Corollary~\ref{T:stabsig}, and a consequence 
for best $n$-term approximations, Corollary~\ref{C:napprox}.
In Section~\ref{s: implementation}, we demonstrate some numerical experiments
that illustrate the stability of ROMP.

\section{Proof of Theorem~\ref{T:stability}}		\label{s: proof}
We shall prove a stronger version of Theorem~\ref{T:stability}, which states that 
{\em at every iteration} of ROMP, either at least $50\%$ of the newly selected coordinates
are from the support of the signal $v$, or the error bound already holds. 

\begin{theorem}[Iteration Invariant of ROMP] \label{T:it}
  Assume $\Phi$ satisfies the Restricted Isometry Condition 
  with parameters $(4n, \e)$ for $\e = 0.01 / \sqrt{\log n}$. 
  Let $v \ne 0$ be an $n$-sparse vector with measurements $x = \Phi v + e$. 
  Then at any iteration of ROMP, after the regularization step where $I$ is the current chosen
  index set, we have $J_0 \cap I = \emptyset$ and (at least) one of the following:
  \renewcommand{\labelenumi}{{\normalfont (\roman{enumi})}}
  \begin{enumerate}
  \item \label{J support}
  $|J_0 \cap \supp(v)| \geq \frac{1}{2}|J_0|$;
  \item \label{error} $\|v|_{\supp(v)\backslash I}\|_2 \leq 100 \sqrt{\log n}\|e\|_2$.
  \end{enumerate} 
  In other words, either at least $50\%$ of the coordinates in the newly selected set $J_0$ 
  belong to the support of $v$ or the bound on the error already holds.
\end{theorem}

We show that the Iteration Invariant implies Theorem~\ref{T:stability} by examining the three possible cases:

{\bf Case 1: (ii) occurs at some iteration.} We first note that since $|I|$ is nondecreasing, if (ii) occurs at some iteration, then it holds for all subsequent iterations. To show that this would then imply Theorem~\ref{T:stability}, 
we observe that by the Restricted Isometry Condition and since $|\supp(\hat{v})| \leq |I| \leq 3n$,
$$
(1-\e)\|\hat{v}-v\|_2 - \|e\|_2 \leq \|\Phi \hat{v} - \Phi v - e\|_2.
$$

Then again by the Restricted Isometry Condition and definition of $\hat{v}$, 
$$
\|\Phi \hat{v} - \Phi v - e\|_2 \leq \|\Phi (v|_I) - \Phi v - e\|_2 \leq (1 + \e)\|v|_{\supp(v)\backslash I}\|_2+\|e\|_2.
$$
Thus we have that
$$
\|\hat{v}-v\|_2 \leq \frac{1+\e}{1-\e}\|v|_{\supp(v)\backslash I}\|_2 + \frac{2}{1-\e}\|e\|_2.
$$
Thus (ii) of the Iteration Invariant would imply Theorem~\ref{T:stability}. 

{\bf Case 2: (i) occurs at every iteration and $J_0$ is always non-empty.} In this case, by (i) and the fact that $J_0$ is always non-empty, the algorithm identifies at least one element of the support in every iteration. Thus if the algorithm runs $n$ iterations or until $|I| \geq 2n$, it must be that $\supp(v) \subset I$, meaning that $v|_{\supp(v)\backslash I} = 0$.  Then by the argument above for Case 1, this implies Theorem~\ref{T:stability}.

{\bf Case 3: (i) occurs at each iteration and $J_0 = \emptyset$ for some iteration.} By the definition of $J_0$, if $J_0 = \emptyset$ then $u = \Phi^* r = 0$ for that iteration. By definition of $r$, this must mean that 
$$
\Phi^*\Phi(v-y) + \Phi^* e = 0.
$$
This combined with Part 1 of Proposition~\ref{P:cons} below (and its proof, see~\cite{NV}) applied with the set $I' = \supp(v)\cup I$ yields
$$
\|v-y+(\Phi^* e)|_{I'}\|_2 \leq 2.03\e\|v-y\|_2.
$$
Then combinining this with Part 2 of the same Proposition, we have
$$
\|v-y\|_2 \leq 1.1\|e\|_2.
$$
Since $v|_{\supp(v)\backslash I} = (v-y)|_{\supp(v)\backslash I}$, this means that the error bound (ii) must hold, so by Case 1 
this implies Theorem~\ref{T:stability}.


\medskip

We now turn to the proof of the Iteration Invariant, Theorem~\ref{T:it}. We will use the following proposition from \cite{NV}.

\begin{proposition}[Consequences of Restricted Isometry Condition \cite{NV}]\label{P:cons} 
  Assume a measurement matrix $\Phi$ satisfies the Restricted Isometry Condition 
  with parameters $(2n, \e)$. Then the following holds.
  \begin{enumerate}
    \item {\em (Local approximation)} 
      For every $n$-sparse vector $v \in \R^d$ 
      and every set $I \subset \{1, \ldots, d\}$, $|I| \le n$, 
      the observation vector $u = \Phi^* \Phi v$ satisfies 
      $$
      \|u|_I - v|_I\|_2 \le 2.03 \e \|v\|_2.
      $$
    \item {\em (Spectral norm)} 
      For any vector $z \in \R^N$
      and every set $I \subset \{1, \ldots, d\}$, $|I| \le 2n$, we have 
      $$
      \|(\Phi^*z)|_I\|_2 \leq (1+\e)\|z\|_2.
      $$
    \item {\em (Almost orthogonality of columns)} 
      Consider two disjoint sets $I,J \subset \{1, \ldots, d\}$, $|I \cup J| \le 2n$.
      Let $P_I, P_J$ denote the orthogonal projections in $\R^N$
      onto $\range(\Phi_I)$ and $\range(\Phi_J)$, respectively. Then 
      $$
      \|P_I P_J\|_{2\rightarrow 2} \leq 2.2 \e.
      $$
  \end{enumerate}
\end{proposition}
The proof of Theorem~\ref{T:it} is by induction on the iteration of ROMP.
The induction claim is that for all previous iterations, the set of newly chosen 
indices is disjoint from the set of previously chosen indices $I$, 
and either (i) or (ii) holds. Clearly if (ii) held in 
a previous iteration, it would hold in all future iterations. Thus we may assume that 
(ii) has not yet held. Since (i) has held at each previous iteration, we must have 
\begin{equation}\label{I2n}|I|\leq 2n.\end{equation}
 
Let $r \ne 0$ be the residual at the start of this iteration, and let $J_0$, $J$ be the sets found by ROMP in this iteration.
As in \cite{NV}, we consider the subspace
$$
H := \range (\Phi_{\supp(v) \cup I})
$$
and its complementary subspaces
$$
F := \range (\Phi_I), \quad 
E_0 := \range (\Phi_{\supp(v) \setminus I}).
$$
The Restricted Isometry Condition in the form of Part~3 of Proposition~\ref{P:cons}
ensures that $F$ and $E_0$ are almost orthogonal. Thus $E_0$ is close to 
the orthogonal complement of $F$ in $H$,
$$
E := F^{\perp}\cap H.
$$

The residual $r$ thus still has a simple description:

\begin{lemma}[Residual]     \label{residual}
  Here and thereafter, let $P_L$ denote the orthogonal projection in $\R^N$ 
  onto a linear subspace $L$. Then
  $$
  r = P_E \Phi v + P_{F^\perp}e.
  $$
\end{lemma}

\begin{proof}
By definition of the residual in the algorithm,
$r = P_{F^\perp} x = P_{F^\perp}(\Phi v + e)$. To complete the proof we need that
$P_{F^\perp}\Phi v = P_E \Phi v$. This follows from the orthogonal decomposition $H=F+E$ and 
the fact that $\Phi v \in H$. 
\end{proof}

Now we consider the signal we seek to identify at the current iteration, 
and its measurements:
\begin{equation}            \label{v0 x0}
  v_0 := v|_{\supp(v) \setminus I}, \quad 
  x_0 := \Phi v_0 \in E_0.
\end{equation}
To guarantee a correct identification of $v_0$, we first state
two approximation lemmas that reflect in two different ways the fact 
that subspaces $E_0$ and $E$ are close to each other.

\begin{lemma}[Approximation of the residual]\label{C:proj}
  We have
  $$
  \|x_0 - r\|_2 \leq 2.2 \e \|x_0\|_2 + \|e\|_2.
  $$
\end{lemma}

\begin{proof}
By definition of $F$, we have 
$\Phi v - x_0 = \Phi(v - v_0) \in F$. 
Therefore, by Lemma~\ref{residual},
$r = P_{E}\Phi v + P_{F^\perp}e = P_{E}x_0 + P_{F^\perp}e$, and so
$$
\|x_0 - r\|_2 = \|x_0 - P_E x_0 - P_{F^\perp}e\|_2 \leq \|P_Fx_0\|_2 + \|e\|_2. 
$$
Now we use Part 3 of Proposition~\ref{P:cons} for the sets $I$ and $\supp(v) \setminus I$
whose union has cardinality at most $3n$ by \eqref{I2n}. It follows that
$$
\|P_Fx_0\|_2 + \|e\|_2 = \|P_F P_{E_0}x_0\|_2 + \|e\|_2 \le 2.2 \e \|x_0\|_2 + \|e\|_2
$$
as desired.
\end{proof}

\begin{lemma}[Approximation of the observation]\label{L:uj}
  Consider the observation vectors 
  $u_0 = \Phi^*x_0$ and $u = \Phi^*r$. Then for any set $T \subset \{1, \ldots, d\}$ with $|T|\leq 3n$,
  $$
  \|(u_0 - u)|_T\|_2 \leq 2.4 \e \|v_0\|_2 + (1+\e)\|e\|_2.
  $$
\end{lemma}

\begin{proof}
Since $x_0 = \Phi v_0$, we have by Lemma~\ref{C:proj} 
and the Restricted Isometry Condition that
$$
\|x_0 - r\|_2 
\le 2.2 \e \|\Phi v_0\|_2 + \|e\|_2 
\le 2.2 \e (1+\e) \|v_0\|_2 + \|e\|_2 
\le 2.3 \e \|v_0\|_2 + \|e\|_2.
$$
To complete the proof, it remains to apply Part 2 of Proposition~\ref{P:cons},
which yields 
$\|(u_0 - u)|_T\|_2 \le (1 + \e)\|x_0 - r\|_2$.
\end{proof}

We next show that the energy (norm) of $u$ when restricted to $J$, and furthermore to 
$J_0$, is not too small. By the regularization step of ROMP, since all selected coefficients have comparable
magnitudes, we will conclude that not only a portion of energy
but also of the {\em support} is selected correctly, or else the error bound must already
be attained. This will be the desired conclusion.

\begin{lemma}[Localizing the energy]\label{C:uj}
  We have $\|u|_J\|_2 \ge 0.8 \|v_0\|_2 - (1+\e)\|e\|_2$.
\end{lemma}

\begin{proof}
Let $S$ = $\supp(v) \setminus I$.
Since $|S| \leq n$, the maximality property of $J$ in the algorithm
implies that 
$$
\|u|_J\|_2 \geq \|u|_S\|_2.
$$
By Lemma~\ref{L:uj}, 
$$
\|u|_S\|_2 \geq \|u_0|_S\|_2 - 2.4\e\|v_0\|_2 - (1+\e)\|e\|_2.
$$
Furthermore, since $v_0|_S = v_0$, by Part 1 of Proposition~\ref{P:cons} we have
$$
\|u_0|_S\|_2 \geq (1 - 2.03\e)\|v_0\|_2.
$$
Putting these three inequalities together, we conclude that
$$
\|u|_J\|_2 \ge (1 - 2.03\e)\|v_0\|_2 - 2.4\e\|v_0\|_2 - (1+\e)\|e\|_2 \ge 0.8 \|v_0\|_2  - (1+\e)\|e\|_2.
$$
This proves the lemma.
\end{proof}

We next bound the norm of $u$ restricted to the smaller set $J_0$, again
using the general property of regularization. 
In our context, Lemma~3.7 of \cite{NV} applied to the vector $u|_J$ yields
$$
\|u|_{J_0}\|_2 \geq \frac{1}{2.5\sqrt{\log n}}\|u|_J\|_2.
$$
Along with Lemma~\ref{C:uj} this directly implies:

\begin{lemma}[Regularizing the energy]\label{C:uj0} 
  We have
  $$
  \|u|_{J_0}\|_2 \ge \frac{1}{4\sqrt{\log n}}\|v_0\|_2 - \frac{\|e\|_2}{2\sqrt{\log n}}.
  $$
\end{lemma}

\medskip

We now finish the proof of Theorem~\ref{T:it}. The claim that $J_0 \cap I = \emptyset$  
follows by the same arguments as in \cite{NV}. 


The nontrivial part of the theorem is its last claim, that either (i) or (ii) holds.
Suppose (i) in the theorem fails. Namely, suppose that 
$|J_0 \cap \supp(v)| < \frac{1}{2}|J_0|$, 
and thus
$$
|J_0 \backslash \supp(v)| > \frac{1}{2}|J_0|.
$$
Set $\Lambda = J_0\backslash\supp(v)$. 
By the comparability property of the coordinates in $J_0$ 
and since $|\Lambda| > \frac{1}{2}|J_0|$, there is a fraction of energy 
in $\Lambda$:
\begin{equation}\label{E:ubig} 
  \|u|_{\Lambda}\|_2 > \frac{1}{\sqrt{5}}\|u|_{J_0}\|_2 
  \ge \frac{1}{4\sqrt{5\log n}}\|v_0\|_2 - \frac{\|e\|_2}{2\sqrt{5\log n}}, 
\end{equation}
where the last inequality holds by Lemma~\ref{C:uj0}.

On the other hand, we can approximate $u$ by $u_0$ as
\begin{equation}                \label{u u0}
  \|u|_{\Lambda}\|_2 
  \le \|u|_{\Lambda} - u_0|_{\Lambda}\|_2 + \|u_0|_{\Lambda}\|_2.
\end{equation}
Since $\L \subset J$, $|J| \leq n$, and using Lemma~\ref{L:uj}, we have
$$
\|u|_{\Lambda} - u_0|_{\Lambda}\|_2 \le 2.4\e\|v_0\|_2 + (1+\e)\|e\|_2.
$$
Furthermore, by definition \eqref{v0 x0} of $v_0$, we have $v_0|_\Lambda = 0$. 
So, by Part 1 of Proposition~\ref{P:cons}, 
$$
\|u_0|_{\Lambda}\|_2 \le 2.03 \e \|v_0\|_2.
$$
Using the last two inequalities and \eqref{u u0}, we conclude that 
$$
\|u|_{\Lambda}\|_2 \le 4.43 \e \|v_0\|_2 + (1+\e)\|e\|_2.
$$
This is a contradiction to~(\ref{E:ubig}) 
so long as 
$$
\e \leq \frac{0.02}{\sqrt{\log n}} - \frac{\|e\|_2}{\|v_0\|_2}.
$$ 
If this is true, then indeed (i) in the theorem must hold. If
it is not true, then by the choice of $\e$, this implies that
$$
\|v_0\|_2 \leq 100\|e\|_2\sqrt{\log n}.
$$
This proves Theorem~\ref{T:it}. Next we turn to the proof of Corollary~\ref{T:stabsig}.
\qed

\section{Approximately sparse vectors and best $n$-term approximations} 
					\label{s: consequences}  

\subsection{Proof of Corollary~\ref{T:stabsig}}


We first partition $v$ so that $x = \Phi v_{2n} + \Phi (v-v_{2n}) + e$. Then since $\Phi$ satisfies the Restricted Isometry Condition with parameters $(8n, \e)$, by Theorem~\ref{T:stability} and the triangle inequality,
\begin{equation}\label{primebnd}
\|v_{2n}-\hat{v}\|_2 \leq 104\sqrt{\log 2n}(\|\Phi (v-v_{2n})\|_2+\|e\|_2),
\end{equation}
The following lemma as in \cite{GSTV2} relates the $2$-norm of a vector's tail to its $1$-norm. An application of this lemma combined with \eqref{primebnd} will prove Corollary~\ref{T:stabsig}.

\begin{lemma}[Comparing the norms]\label{L:ve}
Let $w\in\R^d$, and let $w_m$ be the vector of the $m$ largest coordinates in absolute value from $w$.  Then
$$
\|w-w_m\|_2 \leq \frac{\|w\|_1}{2\sqrt{m}}.
$$
\end{lemma}
\begin{proof}
Let $\mu$ denote the $(n+1)th$ largest entry of $w$. If $\mu=0$ then $w_m=0$ so the claim holds. Thus we may assume this is not the case. Then we have
$$
\frac{\|g-g_m\|_2}{\|g\|_1} \le \frac{\sqrt{\|g-g_m\|_1\|g-g_m\|_\infty}}{\|g_m\|_1+\|g-g_m\|_1} \le \frac{\sqrt{m\mu^2}}{m\mu+m\mu}.
$$ 
Simplifying gives the desired result.
\end{proof}

By Lemma 29 of \cite{GSTV2}, we have 
$$
\|\Phi (v-v_{2n})\|_2 \leq (1+\e)\Big(\|v-v_{2n}\|_2 + \frac{\|v-v_{2n}\|_1}{\sqrt{n}}\Big).
$$
Applying Lemma~\ref{L:ve} to the vector $w=v-v_n$ we then have
$$
\|\Phi (v-v_{2n}\|_2 \leq 1.5(1+\e)\frac{\|v-v_n\|_1}{\sqrt{n}}.
$$
Combined with \eqref{primebnd}, this proves the corollary.

\subsection{Best $n$-term approximation}

We now show that by truncating the reconstructed vector, we obtain a $2n$-sparse vector very close
to the original signal. 
\begin{corollary}\label{C:napprox}Assume a measurement matrix $\Phi$ satisfies the Restricted Isometry Condition 
  with parameters $(8n, \e)$ for $\e = 0.01 / \sqrt{\log n}$. 
  Let $v$ be an arbitrary vector in $\R^d$, let $x=\Phi v + e$ be the measurement vector, and $\hat{v}$ the 
  reconstructed vector output by the ROMP Algorithm. Then
$$
\|v_{2n} - \hat{v}_{2n}\|_2 \leq 477 \sqrt{\log 2n}\Big( \|e\|_2 + \frac{\|v-v_{n}\|_1}{\sqrt{n}}\Big),
$$  
where $z_{m}$ denotes the best $m$-sparse approximation to $z$ 
(i.e. the vector consisting of the largest $m$ coordinates 
in absolute value).
\end{corollary}
\begin{proof}
Let $v_S := v_{2n}$ and $\hat{v_T} := \hat{v}_{2n}$, and let $S$ and $T$ denote the 
supports of $v_S$ and $\hat{v_T}$ respectively.
By Corollary~\ref{T:stabsig}, it suffices to show that 
$\|v_S - \hat{v}_T\|_2 \leq 3\|v_S - \hat{v}\|_2$. 

Applying the triangle inequality, we have
$$
\|v_S - \hat{v}_T\|_2 \leq \|(v_S - \hat{v}_T)|_T\|_2 + \|v_S|_{S\backslash T}\|_2 =: a + b.
$$
We then have
$$
a = \|(v_S - \hat{v}_T)|_T\|_2 \leq \|v_S - \hat{v}_T\|_2
$$
and
$$
b \leq \|\hat{v}|_{S\backslash T}\|_2 + \|(v_S-\hat{v})|_{S\backslash T}\|_2.
$$
Since $|S| = |T|$, we have $|S\backslash T| = |T\backslash S|$. By the definition of $T$, every coordinate of $\hat{v}$ in $T$ is greater than or equal to every coordinate of $\hat{v}$ in $T^c$ in absolute value. Thus we have,
$$
\|\hat{v}|_{S\backslash T}\|_2 \leq \|\hat{v}|_{T\backslash S}\|_2 = \|(v_S - \hat{v})|_{T\backslash S}\|_2.
$$
Thus $b \leq 2\|v_S - \hat{v}\|_2$, and so
$$
a+b \leq 3\|v_S - \hat{v}\|_2.
$$
This completes the proof.
\end{proof}

\medskip

\noindent {\bf Remark.} 
Corollary~\ref{C:napprox} combined with Corollary~\ref{T:stabsig} and \eqref{vbound} implies that we can also estimate a bound on the whole signal $v$:
$$
\|v-\hat{v}_{2n}\|_2 \leq C\sqrt{\log 2n}\Big(\|e\|_2 + \frac{\|v-v_{n}\|_1}{\sqrt{n}}\Big).
$$

\section{Numerical Examples}     \label{s: implementation}

This section describes our experiments that illustrate the stability of ROMP.
We experimentally examine the recovery error using ROMP for both perturbed
measurements and signals. 
The empirical recovery error is actually much better than that given in the theorems. 

First we describe the setup of our experiments. For many values of the ambient dimension $d$, 
the number of measurements $N$, and the sparsity $n$, we reconstruct random signals using ROMP.
For each set of values, we perform $500$ trials. Initially, we generate an $N \times d$ Gaussian measurement matrix $\Phi$. For each trial, independent of the matrix, we generate an $n$-sparse signal $v$ by choosing $n$ components uniformly 
at random and setting them to one.
In the case of perturbed signals, we add to the signal a $d$-dimensional error vector with Gaussian entries. In the case of perturbed measurements, we add an $N$-dimensional error vector with Gaussian entries to the measurement vector $\Phi v$.  We then execute ROMP with the measurement
vector $x = \Phi v$ or $x + e$ in the perturbed measurement case. After ROMP terminates, we output the reconstructed vector $\hat{v}$ obtained from the least squares calculation and calculate its distance from the original signal. 

\begin{figure}[ht] 
  \includegraphics[width=0.8\textwidth,height=3.2in]{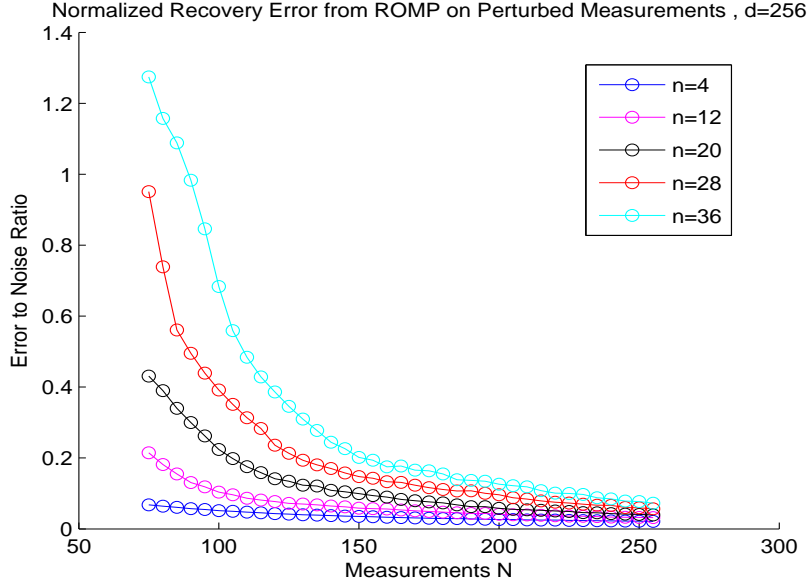}
  \caption{The error to noise ratio $\frac{\|\hat{v} - v\|_2}{\|e\|_2}$ as a function of the number of measurements $N$ in dimension $d=256$ for various levels of sparsity $n$.}\label{fig:meas2}
\end{figure}

Figure~\ref{fig:meas2} depicts the recovery error $\|v - \hat{v}\|_2$ when ROMP was run with perturbed measurements. This plot was generated with $d = 256$ for various levels of sparsity $n$. The horizontal axis represents the number of measurements $N$, and the vertical
axis represents the average normalized recovery error. Figure~\ref{fig:meas2} confirms the results of Theorem~\ref{T:stability}, while also suggesting the bound may be improved by removing the $\sqrt{\log n}$ factor.

\begin{figure}[ht] 
  \includegraphics[width=0.8\textwidth,height=3.2in]{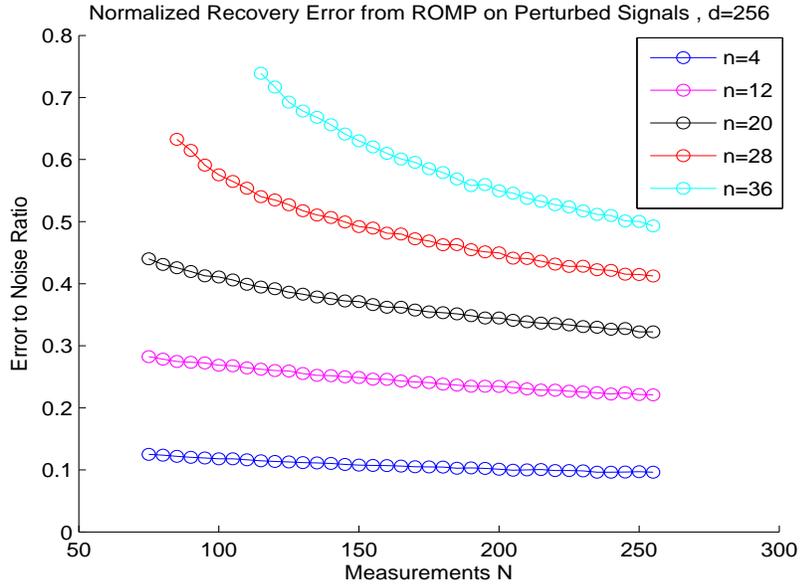}
  \caption{The error to noise ratio $\frac{\|\hat{v} - v_{2n}\|_2}{\|v-v_n\|_1/\sqrt{n}}$using a perturbed signal, as a function of the number of measurements $N$ in dimension $d=256$ for various levels of sparsity $n$.}\label{fig:sig4}
\end{figure}

Figure~\ref{fig:sig4} depicts the normalized recovery error when the signal was perturbed by a Gaussian vector. The figure confirms the results of Corollary~\ref{T:stabsig} while also suggesting again that the logarithmic factor in the corollary is unnecessary. 


\end{document}